\documentclass[reqno]{amsart}
\usepackage{amsmath,amssymb}
\usepackage{graphicx,color}
\usepackage[all]{xy}
\theoremstyle{plain}

\newtheorem{theorem}{Theorem}[section]
\newtheorem{proposition}[theorem]{Proposition}
\newtheorem{lemma}[theorem]{Lemma}
\newtheorem{corollary}[theorem]{Corollary}

\theoremstyle{definition}
\newtheorem{definition}[theorem]{Definition}

\theoremstyle{remark}
\newtheorem{remark}[theorem]{Remark}

\def\E{{\mathcal E}}
\def\B{{\mathcal B}}

\def\M{{\mathcal M}}

\def\Q{{\mathbb Q}}

\def\L{{\mathbb L}}

\def\cat0{\mathrm{cat}_0}

\def\dim{\mathrm{dim}}
\def\ker{\mathrm{ker}}

\begin{document}

\title[]
{On THE GROUP OF  SELF-HOMOTOPY EQUIVALENCES  OF AN $F_0$-SPACE}

\author{Mahmoud Benkhalifa}
\address{Department of Mathematics. College of  Sciences, University of Sharjah.  UAE }

\email{mbenkhalifa@sharjah.ac.ae}



\begin{abstract}
In  \cite{G}, G. Lupton conjectured that the  group of self-homotopy equivalences of an $F_0$-space inducing the identity on  the homotopy  groups is finite.	Thus, the aim of this paper is to establish this conjecture.																								 
\end{abstract}
\maketitle
\section{Introduction}
Let $X$ be a simply connected space with finite dimensional rational homotopy, finite dimensional rational homology (i.e. a rationally elliptic space), and positive Euler characteristic.
The collection of such spaces $X$ is referred to as the class of $F_0$-spaces. Extensively studied
by Halperin, \cite{H}, $F_0$-spaces are rational Poincar\'{e} duality spaces with rational cohomology
$\Q[x_{1},\dots,x_{n}]/(P_{1},\dots,P_{n})$,  where the polynomials $ P_{1},\dots,P_{n} $ form a regular sequence in $\Q[x_{1},\dots,x_{n}]$, i.e., $ P_{1}\neq 0$ and for every $i\geq 1$,  $P_{i}$ is not a zero divisor in  $\Q[x_{1},\dots,x_{n}]/(P_{1},\dots,P_{i-1})$. For instance, 
products of even spheres,  complex  Grassmannian manifolds  and homogeneous spaces $G/H$ such that  rank $G$ = rank $H$ are $F_0$-spaces.

Let $\E(X) $ denote the group  of self-homotopy equivalences of $X$ and let $\E_{\#}(X) $ be its  subgroup  of the elements inducing the identity on the homotopy groups (\cite{B1},\cite{B3}).  

Halperin has conjectured that the rational Serre spectral sequence collapses for any rational
fibration, provided the fiber $X$ is a $F_0$-space. This conjecture, which remains unsolved, can be
rephrased in terms of the (graded Lie algebra of) negative-degree derivations of the rational
cohomology of X (see \cite{M} for more details). Namely:
$$Der_{<0}H^*(X;\Q)=0\Longleftrightarrow \text {  Halperin’s conjecture holds}$$
If we look at the zero-degree derivations of the rational cohomology of $X$, there exists a
correspondence between the decomposable derivations of $Der_{0}H^*(X;\Q)$ and the subgroup  $\E_{\#}(X)$. Hence,
$$Der_{0}H^*(X;\Q)\text { is trivial}\Longrightarrow \E_{\#}(X) \text { is finite}  $$
Motivated by Halperin’s conjecture and this correspondence, Lupton raises the following
question:

\noindent Question(\cite{ak1}, Problem 10): For an  $F_0$-space  $X$, is $\E_{\#}(X)$ finite?

\noindent Thus, the purpose of this paper is to settle this question in the positive using standard tools  of rational homotopy theory which we  refer to \cite{FHT} for a general introduction to these techniques. We recall some of the notation here. By a Sullivan  algebra we mean a free graded commutative algebra $\Lambda V$, for some finite-type graded vector space $V=V^{\geq 2}$, i.e., $\dim\,V^{n}<\infty$ for all $n\geq 2$, together with a differential $\partial$  of degree +1 that is decomposable, i.e., satisfies $\partial: V \to \Lambda^{\geq 2} V$.  Here $\Lambda^{\geq 2} V$ denotes the graded vector space spanned by all  the monomials $v_{1}\dots v_{r}$ such that $v_{1},\dots, v_{r}\in V$ and $r\geq 2$.

Every simply connected   space $X$ with rational cohomology of finite-type has a corresponding  Sullivan  algebra  called the Sullivan model of $X$,   unique up to isomorphism, that encodes the rational  homotopy type of  $X$. In particular we have  
$$ V^*\cong \text{Hom}(\pi_{*}(X)\otimes \Bbb Q,\Bbb Q),\,\,\,\,\,\,\,\,\,\,\,\,\,\,\,\,\,\,\,\,\,\,\,\,\,H^{*}(\Lambda V)\cong H^{*}(X,\Bbb Q).$$

By a free differential graded Lie algebra $(\L(W),\delta)$ (DGL for short),  we mean a free graded Lie algebra  $\L(W)$, for some finite-type  vector space $W=(W_{\geq 1})$, together with a decomposable differential $\delta$  of degree -1, i.e., $\delta (W) \subset \L^{\geq 2} (W)$.  Here $\L^{\geq 2} (W)$ denotes the graded vector space spanned by all  the brackets of lengths  $\geq 2$.

Every  simply connected   space $X$ with rational cohomology of finite-type has  a corresponding DGL $(\L(W),\delta)$ called the Quillen model of $X$,  unique up to isomorphism, which determines completely the rational homotopy type of $X$.  In particular we have  
\begin{equation*}\label{27}
	  W_*\cong H_{*+1}(X;\Bbb Q),\,\,\,\,\,\,\,\,\,\,\,\,\,\,\,\,\,\,\,H_{*}(\L(W))\cong \pi_{*+1}(X)\otimes \Bbb Q.
\end{equation*}

This work consists of five sections, the first one being the introduction. Section 2
is devoted to state some results on the  notion  of DGL-homotopy  as well as the properties of an $F_0$-space $X$ notably, if   $(\L(W),\delta)$ is its  Quillen model, then  we introduce the group $\E_{\#}(\L(W))$ of the self-homotopy equivalences  of  $(\L(W),\delta)$ constituting on the elements $[\alpha]$ satisfying $H_{*}(\alpha)=\rm{id}$.  In Section 3, we prove that if  $[\alpha]\in\E_{\#}(\L(W))$, then $\alpha$ is homotopic to DGL-map $\tilde{\alpha}$ satisfying  $\tilde{\alpha}(W)=W$. In Section 4 and 5, we focus on studying  the properties of $(\L(W),\delta)$ to show  that  $\E_{\#}(\L(W))$ is trivial. Consequently,  by virtue of  the   localization theorem of Maruyama \cite{Mar}, we  derive that $\E_{\#}(X)$ is finite.
\section{Preliminaries}
\subsection{Homotopy between DGL-maps (see \cite[\S 21]{FHT})}
	Let $(\L(W),\delta)$ be a  DGL. Define the DGL $\L(W, sW, W'),D)$ with $W\cong W'$ and $(sW)_{i}=W_{i-1}$. The differential $D$ is given by
\begin{equation}\label{g3}
D(w)=\partial(w),\,\,\,\,\,\,\,\,\,\,\,\,\,\,\,\,\,\,\,\,\,\,  D(sw)=w',\,\,\,\,\,\,\,\,\,\,\,\,\,\,\,\,\,\,\,  D(w')=0.
\end{equation}
Define  $S$ as the derivation of degree +1 on $\L(W, sW, W')$ given by $$S(w)=sw,\,\,\,\,\,\,\,\,\,\,\,\,\,\,\,\,\,\,\,\,\,\,S(sw)=S(w')=0.$$
A  homotopy   between two DGL-maps  $\alpha,\alpha':(\L(W),\delta)\to (\L(W),\delta)$  is DGL-map
$$F \colon   (\L(W, sW, W'),D)\to (\L(W),\delta),$$
such that $F(w)=\alpha(w)$ and $F\circ e^{\theta}(w)=\alpha'(w)$,  where
$$e^{\theta}(w)=w+w'+\underset{n\geq 1}{\sum} \frac{1}{n!}(S\circ D)^{n}(w), \text{\,\,\,\,and \,\,\,\,} \theta=D\circ S+S\circ D.$$
Thus, the notion of DGL-homotopy allows us to define  the following group. 
\begin{definition}\label{d3}
Let $\E_{\#}(\L(W))$ denote the group of self-homotopy equivalences  of  $(\L(W),\delta)$ constituting with the elements $[\alpha]$ satisfying $H_{*}(\alpha)=\rm{id}$, where $$H_{*}(\alpha):H_{*}(\L(W_{\rm{}}))\to H_{*}(\L(W_{\rm{}})).$$ 
\end{definition}
By virtue of the  properties of the   model of Quillen and  the   localization theorem of Maruyama \cite{Mar},  we deduce that if $X$ is an $F_0$-space, then we have 
\begin{equation}\label{f13}
	\E_{\#}(X)\otimes \Q\cong\E_{\#}(\L(W)).
\end{equation}
Thus, the group $\E_{\#}(X)$ is finite if and only if the group $\E_{\#}(\L(W))$ is trivial.

\medskip
 Later on we will need the following two lemmas.
\begin{lemma}
	\label{l00}   Let   $\alpha,\tilde{\alpha}\colon   (\L(W_{\leq n}),\delta) \to (\L( W_{\leq n}),\delta)$  be two DGL-maps such that 
	$$\alpha(w)=\tilde{\alpha}(w)+y \,\text{ on }   W_{n} \,\,\,\,
	\text{ and }\,\,\,\, \alpha=\tilde{\alpha}
	\,\,\text{ on }  W_{ \leq n-1}.$$
	Assume that $y=\delta(z)$,
	where $z\in \L(W_{\leq n})$. Then $\alpha$ and $\tilde{\alpha} $
	are  homotopic.
\end{lemma}
\begin{proof}
Define $F \colon   (\L(W_{\leq n}, sW_{\leq n}, W'_{\leq n}),D)\to (\L(W_{\leq n}),\delta)$ by setting
	\begin{eqnarray}
		\label{f6}
		F(w) \hspace{-2mm}&=&\hspace{-2mm} \alpha(w), \,\,\,\,\,F(w')=-y \hbox{\ \ and \ \ }F(sw)=-z \hbox{\ \ for \ \ } w \in W_{n}, \\
		F(w)\hspace{-2mm}&=&\hspace{-2mm} \alpha(w), \,\,\,\,\,\,F(w')=0 \,\,\,\,\hbox{\ \ and \ \ }F(sw)=0\,\,\,\, \hbox{\ \ for \ \ } w \in   W_{\leq n-1}\nonumber.
	\end{eqnarray}
Let $w \in W_{n}$, by considering the relations (\ref{g3}), (\ref{f6}) and as $\delta(w)\in\L(W_{\leq n-1})$, we get 
	$$\delta F(w)=\delta \alpha(w),\,\,\,\,\,\,\,\,\,\,\,\,\,\,\,FD(w)=F(\delta(w))=\alpha\delta(w).$$
Moreover,  a straightforward computation shows 
\begin{eqnarray}
	\label{f7}
\delta F(w') \hspace{-2mm}&=&\hspace{-2mm}\delta(-y)=-\delta(\delta(z))=0,\,\,\,\,\,\,\,\,\,\,\,\,\,\,\,FD(w')=F(0)=0,\nonumber\\
	\delta F(sw) \hspace{-2mm}&=&\hspace{-2mm}\delta(-z)=-y,\,\,\,\,\,\,\,\,\,\,\,\,\,\,\,\,\,\,\,\,\,\,\,\,\,\,\,\,\,\,\,\,\,\,\,\,\,\,\,\,\,FD(sw)=F(w')=-y,\nonumber
	\end{eqnarray}
implying that $F$ is a DGL-map. Next, on the one hand,  from (\ref{f6}), we have $F(w)=\alpha(w)$ for every $w\in W$.   On the other hand, by  expanding the expression $(S\circ\partial)^{n}(w)$ leads to linear combinations of brackets involving the generators $sw$, where $w\in W_{\leq n-1}$. Since in this case  $F(sw)=0$,  it follows that 
$\underset{n\geq 1}{\sum} \frac{1}{n!}F(S\circ D)^{n}(w)=0$. Consequently, 	 we obtain
\begin{eqnarray}
\label{f11}
F\circ e^{\theta}(w)\hspace{-2mm}&=&\hspace{-2mm}F(w)+F(w')= \alpha(w)-y=\tilde{\alpha}(w)\,\,,\,\,\,\,\,\,\,\text{ 	if $w\in W_{n}$ },\nonumber\\
F\circ e^{\theta}(w)\hspace{-2mm}&=&\hspace{-2mm}F(w)+F(w')= \alpha(w),\,\,\,\,\,\,\,\,\,\,\,\,\,\,\,\,\,\,\,\,\,\,\,\,\,\,\,\,\,\,\,\,\,\,\,\,\,\,\,\,\text{ 	if $w\in W_{\leq n-1}$ }.\nonumber
\end{eqnarray}	
But by hypothesis we have $\alpha(w)=\tilde{\alpha}(w)$ on $W_{\leq n-1}$, so  	for all  $w\in W$ we have   $F\circ e^{\theta}(w)=\tilde{\alpha}(w)$ implying that $F$ is the needed  homotopy.
\end{proof} 
\begin{lemma}
	\label{c1} Let   $\alpha,\beta \colon   (\L(W_{\leq n}),\delta) \to (\L(W_{\leq n}),\delta)$  be two DGL-maps such that 
	\begin{eqnarray}\label{G1}
		\alpha(w)&=&\beta(w)+y,\,\,\,\,\, w\in W_{n},\,\,\,\,\,\, y\in \L_{n}(W_{\leq n-1}),\nonumber\\
		\alpha&\simeq&\beta,\,\,\,\,\,\,\,\,\,\,\,\,\,\mathrm{ 	on }\,\, \,\L_{}(W_{\leq n-1})\nonumber.
	\end{eqnarray}
	There is a cycle $y'\in \L_{n}(W_{\leq n-1})$ such that  $\alpha$ is homotopic to the following DGL-map
	\begin{eqnarray}\label{g2}
		\alpha'(w)&=&\beta(w)+y',\,\,\,\,\, w\in W_{n},\nonumber\\
		\alpha'&=&\beta,\,\,\,\,\,\,\,\,\,\,\,\,\,\mathrm{ 	on }\,\, \,\L_{}(W_{\leq n-1}).
	\end{eqnarray}
\end{lemma}
\begin{proof}
	 Since $\alpha$ and $\beta$ are homotopic on $\L(W_{\leq n-1})$, there exits a homotopy  $$F \colon   (\L(W_{\leq n-1}, sW_{\leq n-1}, W'_{\leq n-1}),D)\to (\L(W_{\leq n-1}),\delta),$$
	such that
	\begin{equation}\label{g1}
		F(w)=\beta(w),\,\,\,\,\,\,\,\,\,\,\,F\circ e^{\theta}(w)=\alpha(w),\,\,\,\,\,\,\,\,\,\,\,\forall w\in W_{\leq n-1}.
	\end{equation}
	Therefore for $w\in W_{n}$, the element $F\Big(\underset{n\geq 1}{\sum} \frac{1}{n!}(S\circ D)^{n}(w)\Big)$ is a well-defined element in $\L_{n}(W_{\leq n-1})$. Thus we define
	\begin{equation}\label{3}
		y'=y-F\Big(\underset{n\geq 1}{\sum} \frac{1}{n!}(S\circ D)^{n}(w)\Big).
	\end{equation}
	Now, by hypothesis we have
	\begin{equation}\label{4}
\delta(\beta(w))+\delta(y)=	\delta\alpha(w)=\alpha_{}(\delta(w))=F\circ e^{\theta}(\delta(w))=F\circ e^{\theta}(D(w)).
	\end{equation}
	But $e^{\theta}$ is a DGL-automorphism, so
	\begin{eqnarray}\label{5}
		F\circ e^{\theta}(D(w))&=&F\circ D (e^{\theta}(w))=F\circ D\Big(w+w'+\underset{n\geq 1}{\sum} \frac{1}{n!}(S\circ D)^{n}(w)\Big)\nonumber\\
		&=&F(D(w))+F(D(w'))+F\circ D\Big(\underset{n\geq 1}{\sum}\frac{1}{n!}(S\circ D)^{n}(w)\Big)\nonumber\\
		&=&F(\delta(w))+\delta\circ F\Big(\underset{n\geq 1}{\sum}\frac{1}{n!}(S\circ D)^{n}(w)\Big)\nonumber\\
		&=&\beta(\delta(w))+\delta\Big(\underset{n\geq 1}{\sum}\frac{1}{n!}F(S\circ D)^{n}(w)\Big)\nonumber\\
		&=&\delta(\beta(w))+\delta\Big(\underset{n\geq 1}{\sum}\frac{1}{n!}F(S\circ D)^{n}(w)\Big).
	\end{eqnarray}
Here we use the facts that $D(w')=0$ by (\ref{g3}), $F\circ D=\delta\circ F$ and $F(\delta(w))=\beta(\delta(w))$ because 
$\delta(w)\in \L(W_{\leq n-1})$ and $F=\beta$  on $W_{\leq n-1}$ by (\ref{g1}). Comparing (\ref{4}) and (\ref{5}) we get
	\begin{equation*}\label{6}
		\delta(y)=\delta\Big(\underset{n\geq 1}{\sum}\frac{1}{n!}F(S\circ D)^{n}(w)\Big),
	\end{equation*}
	which implies according to (\ref{3}) that $\delta(y')=0.$	
\medskip

 Now define  $G \colon   (\L(W_{\leq n}, sW_{\leq n}, W'_{\leq n}),D)\to (\L(W_{\leq n}),\delta)$  by setting
	\begin{eqnarray}
		\label{f16}
		G(w) \hspace{-2mm}&=&\hspace{-2mm} \alpha'(w),\,\,\,\,\,\,\,\,\,\,\,\,\,\,\,\,\,G(w')=G(sw)=0,\,\,\,\,\,\,\,\,\,\,\,\,\,\,\,\,\,\,\hbox{ for  } w \in W_{n}, \nonumber\\
		G\hspace{-2mm}&=&\hspace{-2mm} F, \,\,\,\,\,\,\,\,\,\,\,\,\,\,\,\,\,\,\,\,\,\,\,\,\, \hbox{ on   }   W_{\leq n-1}\nonumber.
	\end{eqnarray}
Let us consider the relations  (\ref{g3}). A simple computation shows  that 
	$$\delta(G(w))=\delta(\alpha'(w)),\,\,\,\,\,\,\,\,\,\,\,\,\,\,\,\,\,\,\,\,\,\,\,\, G(D(w))=G(\delta(w)).$$
As  $\delta(w)\in \L_{n}(W_{\leq n-1})$,  it follows that $G(\delta(w))=F(\delta(w))$ and by (\ref{g2}),   (\ref{g1}) we get $F(\delta(w))=\beta(\delta(w))=\alpha'(\delta(w))$. As a result $\delta(G(w))=G(D(w))$. Also by taking into consideration the relations (\ref{g3}), we obtain 
$$\delta(G(w'))=GD(w')=0,\,\,\,\,\,\,\,\,\,\,\,\,\,\,\,\,\,\delta(G(sw))=0,\,\,\,\,\,\,\,\,\,\,\,\,\,\, GD(sw)=G(w')=0,$$	
proving that $G$ is a DGL-map satisfying $G(w)=\alpha'(w)$ for all $w\in W_{\leq n}$. Moreover, we have
	$$G\circ e^{\theta}(w)=G\Big(w+w'+\underset{n\geq 1}{\sum} \frac{1}{n!}(S\circ D)^{n}(w)\Big)=G(w)+G\Big(\underset{n\geq 1}{\sum}\frac{1}{n!}(S\circ D)^{n}(w)\Big).$$
	As $\underset{n\geq 1}{\sum} \frac{1}{n!}(S\circ D)^{n}(w)\in \L_{n}(W_{\leq n-1})$ and $F= G$ on $W_{\leq n-1}$, it follows that
	$$G\circ e^{\theta}(w)=\alpha'(w)+F\Big(\underset{n\geq 1}{\sum}\frac{1}{n!}(S\circ D)^{n}(w)\Big)=(\beta(w)+y')+(y-y')=\alpha(w).$$
	Here we use  (\ref{3}). Consequently, $\alpha$ and $\alpha'$ are homotopic.
\end{proof}
\subsection{Whitehead exact sequence of a DGL} 
Let    	$(\L(W),\delta)$ be a DGL. If 
$$j_n:H_{n}(\L (W_{\leq n}))\to W_n,\,\,\,\,\,\,\,\,\,\,\,\,\,\,\,\,\,\,\,\,\,\,\,\,\,\,j_n(\{w+y\})=w,$$ 
where $w\in W_{n}$,   $y\in \L_{n} (W_{\leq n-1})$ and where $\{w+y\}$ denote the homology class of the cycle $w+y$, then we define  the graded vector space $\Gamma_{*}$ by setting 
\begin{equation}\label{a1}
	\Gamma_{n}=\ker \big(H_{n}(\L (W_{\leq n}))\overset{j_n}{\longrightarrow} W_n\big),\,\,\,\,\,\,\,\,\,\,\,\,\,\,\,\,\,\,\,\,\,\,\,\,\,\forall n\geq 2.
\end{equation}
To  every  DGL	$(\L(W),\delta)$, we can assign (see  \cite{B3,B1} for more details) the following long  exact   sequence
\begin{equation}\label{23}
	\cdots \rightarrow W_{n+1}\overset{b_{n+1}}{\longrightarrow }%
	\Gamma_{n}\rightarrow H_{n}(\L(W))\rightarrow
	W_{n}
	\overset{b_{n}}{\longrightarrow }\cdots
\end{equation}
called the Whitehead exact sequence of $(\L(W),\delta)$. Here  $b_n(w)=\{\delta(w)\}$, where $\{\delta(w)\}$  denotes
the homology class of $\delta(w)$ in 
$\L_{n-1} (W_{\leq n-1}).$

\subsection{Elliptic spaces}
Recall that 
a  simply connected  space $X$ is called rationally elliptic if it satisfies $\dim\,(\pi_{*}(X)\otimes \Bbb Q)<\infty$ and $\dim\,H^{*}(X,\Q)<\infty$ (\cite{FHT}, \S32). The following result mentions some important properties of rationally  elliptic spaces.
\begin{proposition}\label{t13}(\cite{FHT} Proposition  \rm{32.6} and \rm{32.10}).
	 If $(\L(W), \delta)$ is the Quillen model of a rationally   elliptic space of formal dimension $M$, then 
	 \begin{itemize}
	 	\item $\dim\,W_{M-1}=1$ and $W_{i}=0$ for all $i\geq M.$
	 	\item $\sum_{i\geq 1}^{}(2i+1)\,\dim\,H_{2i}(\L(W))-\sum_{i\geq 1}^{}(2i)\,(\dim\,H_{2i-1}(\L(W))-1)= M.$
	 \end{itemize}
Furthermore, the following   statements are equivalent 
	\begin{enumerate}
		\item $X$  is  an $F_0$-space.
		\item $\dim\,H_{\rm{even}}(\L(W))=\dim\,H_{\rm{even}}(\L(W)). $
		\item $W_{\rm{even}}=0.$
	\end{enumerate}
\end{proposition}
\begin{remark}\label{r4}
According to  Proposition \ref{t13}, 	the formal dimension  of  an $F_0$-space  must be   an  even integer.
\end{remark}
\section{Properties of  the group    $\E_{\#}(\L(W))$}
The purpose of this section is to study the properties of  the group    $\E_{\#}(\L(W))$, introduced  in definition  \ref{d3},  in the case where the DGL $(\L(W_{\rm{}}),\delta)$ is the Quillen  model of an $F_0$-space. 
 
As it is stated in the introduction,  an $F_0$-space   is an elliptic space such that its  rational cohomology  is a graded algebra on the form $\Q[x_{1},\dots,x_{n}]/(P_{1},\dots,P_{n})$, where  the polynomials $ P_{1},\dots,P_{n} $ form a regular sequence in $\Q[x_{1},\dots,x_{n}]$. 

\noindent In \cite{H}, it is shown that the  Sullivan model of  an  $F_0$-space  is  given by
\begin{equation*}\label{222}
	(\Lambda V,\partial)=(\Lambda( x_{1},\dots,x_{n};y_{1},\dots,y_{n}),\partial)\,\,,\,\,\partial(x_{i})=0\,\,,\,\,\partial(y_{i})=P_{i}\,\,,\,\,1\leq i\leq n,
\end{equation*}
where  the generator $x_{1},\dots,x_{n}$ are  of even degrees and   $y_{1},\dots,y_{n}$ are of  odd degrees.
 
\medskip
It well-known that  $F_0$-spaces  are formal  (see \cite{H}, theorem 5), i.e.,   there exists a quasi-isomorphism $\M(X)\to (H_{*}(X,\Bbb Q),0)$. As a result,   the differential    of the  Quillen model $(\L(W),\delta)$ is purely quadratic, i.e., $\delta(W_{})\subset [W,W]$  (see \cite{N}, proposition 3.2).   Moreover, taking into account that $W_{\rm{even}}=0$,  we deduce that $W=W_{\rm{odd}}$.
 
\medskip
\begin{remark}\label{r0}
Recall that we have $V^{\rm{even}}\cong H_{\rm{odd}}(\L(W_{\rm{}}))$, therefore,  to each  $x_{i}\in V^{\rm{even}}$  corresponds a homology class $\{w_i+q_i\}\in H_{\rm{odd}}(\L(W_{\rm{}}))$ such that $w_i$ is indecomposable and    $q_i$ is decomposable.  Since   $\delta(w_{i})=-\delta(q_{i})$, it follows that $\delta(w_{i})$ has bracket length greater or equal  than 3. But   $\delta$ is purely quadratic, it follows that $q_i=0$. As a result,  $H_{\rm{odd}}(\L(W_{\rm{}}))$ is generated by  $w_1,\dots,w_n$.
\end{remark}
 \begin{proposition}\label{pp1}
Let $(\L(W_{\rm{}}),\delta)$ be the Quillen  model of an $F_0$-space $X$. 	Then the graded vector space    $\Gamma_{\rm{odd}}$, defined in \rm{(\ref{a1})}, is trivial.
\end{proposition}
\begin{proof}
Assume there is  $\{z\}\neq 0\in \Gamma_{\rm{odd}}$. Since $W_{\rm{even}}=0$,   the exact sequence (\ref{23})implies that $\{z\}\in H_{\rm{odd}}(\L(W))$ which is impossible as $z$ is decomposable due to Remark \ref{r0}.
\end{proof}                                                      
Let us consider the Quillen  model $(\L(W_{\rm{}}),\delta)$ of an $F_0$-space $X$ of formal dimension $M$. By virtue of   Proposition \ref{t13}, we can write 
\begin{equation*}\label{s24}
W_{}=W_{r_1}\oplus\dots\oplus W_{r_m}\oplus W_{M-1},\,\,\,\,\,\,\,\,\,\,\,\,\,\,\,\,\,\,\,\,\,\,\,\,\,\,\,\,\,r_1<\dots <r_m< M-1, 
\end{equation*}
\begin{equation}\label{s23}
	W_{r_i}=\big\langle w_{(1,r_i)},\dots,w_{(n_i,r_{i})}\big\rangle,\,\,\,\,\,\,\,\,\,\,\,\,\,\,\,\,1\leq i\leq m,\,\,\,\,\,\,\,\,\,\,\,\,\,\,\,\,\,\,\,W_{M-1}=\langle \mu\rangle. 
\end{equation}
If $[\alpha]\in\E(\L(W))$, then   for every $1\leq i\leq m$ and $1\leq j\leq n_i$, let us write
\begin{eqnarray}\label{f5}
	\alpha(w_{(j,r_i)})&=&\sum_{s_i=1}^{n_i}\lambda_{(j,r_i),s_i}w_{(j,r_i)}+A_{(j,r_i)},\,\,\,\,\,\,\,\,\,\,\,\,\,\,\,\,\,\,\,\,A_{(j,r_i)}\in \L^{\geq 3 }(W_{\leq
		 r_{i-1}}),\nonumber\\
		\alpha(\mu)&=&a\mu+A_{\mu},\,\,\,\,\,\,\,\,\,\,\,\,\,\,\,\,\,\,\,\,\,\,\,\,\,\,\,\,\,\,\,\,\,\,\,\,\,\,\,\,\,\,\,\,\,\,A_{\mu}\in \L^{\geq 3 }(W_{\leq
		M-2}),
\end{eqnarray}
where all the coefficients  $\lambda_{(r_i,j),s_i},a$ are rationals.
\medskip

	Set  $\tilde{\alpha}(w_{(j,r_i)})=\sum_{s_i=1}^{n_i}\lambda_{(j,r_i),s_i}w_{(j,r_i)}$,  then  (\ref{f5}) becomes $$\alpha(w_{(j,r_i)})=\tilde{\alpha}(w_{(j,r_i)})+A_{(j,r_i)}.$$
	 Note that $\tilde{\alpha}(w_{(j,r_i)})\in W_{r_i}$.  Moreover,  if $l(A_{(j,r_i)})$ denotes the bracket length of $A_{(j,r_i)}$, then  $l(A_{(j,r_i)})\geq 3$ because $\vert
	 A_{(j,r_i)}\vert$ is odd and $W=W_{\rm{odd}}$.   

\begin{theorem}\label{p2}
	Let  $X$ be an  $F_0$-space and let  $(\L(W_{\rm{}}),\delta)$ be its Quillen model. If    $[\alpha]\in\E_{\#}(\L(W))$, then  $\alpha$ is homotopic to the  DGL-map 
	$\tilde{\alpha}.$  Here $\E_{\#}(\L(W))$ is  defined in \rm(\ref{d3}).
\end{theorem}
\begin{proof}
Let  $[\alpha]\in\E_{\#}(\L(W))$ and $\alpha_{r_k}:(\L(W_{\leq r_k}),\delta)\to (\L(W_{\leq r_k}),\delta)$, the restriction of $\alpha$ to $\L(W_{\leq r_k})$.  Since $H_{r_1}(\alpha)=id_{H_{r_1}(\L(W_{\rm{odd}}))}=id_{W_{r_1}}$, 
we deduce that $\alpha_{r_3}=id$ on $W_{r_1}$. 

First, from the relation  (\ref{f5}),  we have
	$$\alpha_{r_2}(w_{(j,r_2)})=\tilde{\alpha}_{r_2}(w_{(j,r_2)})+A_{(j,r_2)},\,\,\,\,\,\,l(A_{(j,r_2)})\geq 3,\,\,\,\,\,\,\,\,\,\,\,\,\,
	\alpha_{r_2}=id\,\,,\,\,\mathrm{ 	on }\,\,\L_{}(W_{ r_1}),$$
implying that
$$\delta\alpha_{r_2} (w_{(j,r_2)})=\delta(\tilde{\alpha}_{r_2}(w_{(j,r_2)}))+\delta(A_{(j,r_2)}).$$
Next, as 	$\delta(t_{(j,r_2)})\in\L_{}(W_{r_1})$,  we get 
$$\alpha_{r_2}\delta(w_{(j,r_2)})=\delta(w_{(j,r_2)}).$$
Since $\delta\alpha_{r_2} =\alpha_{r_2}\delta$,  $l(A_{(j,r_2)})\geq 3$ and  $\delta$ is purely quadratic, it follows that $\delta(A_{(j,r_2)})=0$ for every  $1\leq j\leq n_2$.  As a result,    the homology class $\{A_{(j,r_2)}\}$ belongs to $ \Gamma_{r_2}$ which is, by proposition   \ref{pp1}, trivial as $r_2$ is odd, therefore  $A_{(j,r_2)}$ is a boundary. Now applying lemma \ref{l00},  it follows that $\alpha_{r_2}$ and $\tilde{\alpha}_{r_2}$ are homotopic on $\L(W_{\leq r_2}).$

\noindent Assume by induction that $\alpha_{r_{k-1}}$ and $\tilde{\alpha}_{r_{k-1}}$ are homotopic on $\L(W_{\leq {r_{k-1}}})$. Therefore using (\ref{f5}) we get
\begin{eqnarray}\label{g10}
	\alpha_{r_k}(w_{(j,r_k)})&=&\tilde{\alpha}_{r_k}(w_{(j,r_k)})+A_{(j,r_k)},\,\,\,\,\,\,\,\,\,\,\,\,\,\,\,\,l(A_{(j,r_k)})\geq 3,\nonumber\\
\alpha_{r_{k-1}}&\simeq&	\tilde{\alpha}_{r_{k-1}},\,\,\,\,\,\,\,\,\,\,\,\,\,\,\,\,\,\,\,\,\,\,\,\,\,\,\,\,\,\,\,\,\,\,\,\,\,\,\,\,\,\,\,\,\,\mathrm{ on }\,\, \,\L_{}(W_{\leq r_{k-1}}).\nonumber
\end{eqnarray}
Due to lemma \ref{c1}, we deduce that there is a cycle $A'_{(j,r_k)}$ such that $l(A'_{(j,r_k)})\geq 3$ and $\alpha_{r_k}$ is homotopic to the DGL-map $\alpha'_{r_k}$ given  by
\begin{eqnarray}\label{g100}
	\alpha'_{r_k}(w_{(j,r_k)})&=&\tilde{\alpha}_{r_k}(w_{(j,r_k)})+A'_{(j,r_k)},\,\,\,\,\,\,\,\,\,\,\,\,\,\,\,\,\,\,\,\,\,\,\,\,l(A_{(j,r_k)})\geq 3,\nonumber\\
	\alpha'_{r_{k-1}}&=&	\tilde{\alpha}_{r_{k-1}},\,\,\,\,\,\,\,\,\,\,\,\,\,\,\,\,\,\,\,\,\,\,\,\,\,\,\,\,\,\,\,\,\,\,\,\,\,\,\,\,\,\,\,\,\,\,\,\,\,\,\,\,\,\,\,\,\,\,\,\,\,\,\,\mathrm{ on }\,\, \,\L_{}(W_{\leq r_{k-1}}).\nonumber
\end{eqnarray}
The cycle $A'_{(j,r_k)}$ defines a homology class $\{A'_{(j,r_k)}\}$ belonging to $ \Gamma_{\rm{odd}}$ which is trivial by  \ref{pp1} because $\vert A'_{(j,r_k)}\vert=r_k=\rm{odd}$. Therefore, from lemma \ref{l00},  we deduce that $\alpha'_{k}\simeq\tilde{\alpha}_{r_k}$  and so are   $\alpha_{r_k}$ and $\tilde{\alpha}_{r_k}$. Hence,  $\alpha_{}\simeq\tilde{\alpha}_{}$.
\end{proof}
As a consequence of Theorem \ref{p2}, we deduce the following fact
\begin{corollary}\label{c3}
Let  $X$ be an  $F_0$-space and let  $(\L(W_{\rm{}}),\delta)$ be its Quillen model. If $[\alpha]\in\E_{\#}(\L(W))$, then   for every $1\leq s\leq m$ we have
	$\alpha(W_{r_s})=W_{r_s}$ and $\alpha(\mu)=a\mu$, where  $a$ is a non-zero   rational. 
\end{corollary}
\begin{proof}
It follows from Theorem \ref{p2} and the relations 	(\ref{f5}).
\end{proof}
\begin{corollary}\label{c4}
Let  $X$ be an  $F_0$-space and let  $(\L(W_{\rm{}}),\delta)$ be its Quillen model. 	If $[\alpha]\in\E_{\#}(\L(W))$, then   for every indecomposable cycle $w_{(j,r_s)}\in W$, we have $\alpha(w_{(j,r_s)})=w_{(j,r_s)}.$
\end{corollary}
\begin{proof}
By virtue of (\ref{f13}), if $[\alpha]\in\E_{\#}(\L(W))$, then  $H_{*}(\alpha)=id_{H_{*}(\L(W))}$. Therefore, since  $w_{(j,r_s)}$ is a cycle we get
$$H_{*}(\alpha)(\{w_{(j,r_s)}\})=\{w_{(j,r_s)}\},$$
implying $\alpha(w_{(j,r_s)})-w_{(j,r_s)}$ is a boundary in $(\L(W),\delta)$. As $\delta$ is purely quadratic, it follows that  $\alpha(w_{(j,r_s)})=w_{(j,r_s)}.$
\end{proof}
\section{Properties of the Quillen model of an  $F_0$-space}
Let  $X$ be an  $F_0$-space of formal dimension $M$ and let 
\begin{equation*}\label{222}
	(\Lambda V,\partial)=(\Lambda( x_{1},\dots,x_{n};y_{1},\dots,y_{n},\partial)\,\,,\,\,\partial(x_{i})=0\,\,\,,\,\,1\leq i\leq n,
\end{equation*}
 be its  Sullivan model and $(\L(W),\delta)$  its Quillen model. Assume  that
 \begin{equation*}\label{x37}
 \vert x_{1}\vert\leq\dots\leq\vert x_{n}\vert.
 \end{equation*}
Recall that a basis   of $W_{r_s}$ is  given by (see (\ref{s23}))  
\begin{equation*}\label{x14}
	W_{r_s}=\big\langle w_{(1,r_s)},\dots,w_{(n_s,r_{s})}\big\rangle,\,\,\,\,\,\,\,\,\,\,\,\,\,\,\,\,1\leq s\leq m,\,\,\,\,\,\,\,\,\,\,\,\,\,\,\,\,\,\,\,W_{M-1}=\langle \mu\rangle. 
\end{equation*}
To each generator $w_{(j,r_s)}$ corresponds a non-trivial cohomology class    
$\{x^{i_1}_{1}\dots x^{i_n}_{n}\}$  such that 
\begin{equation}\label{x12}
  r_s=i_1\vert x_{1}\vert+\dots+i_n\vert x_{n}\vert -1,\,\,\,\,\,\,\,\,\,\,\,\,\,\,\,\,\,\, i_1\geq 0, \dots,i_n\geq 0.
\end{equation}
The differential is given by
\begin{equation}\label{y2}
	\delta(w_{(j,r_s)})=\sum_{}^{} \lambda_{(i,t)}[w_{(i,r_p)},w_{(t,r_q)}],\,\,\,\,\,\,\,\,\,\,\,\,\, r_p\leq r_q\,\,,\,\,r_p+r_q=r_s-1,
\end{equation}
where  $\lambda_{(i,t)}\in \Q$ and where the generators    $w_{(i,r_p)}\in W_{r_p}$ and $w_{(t,r_q)}\in W_{r_q}$  correspond respectively to the non-trivial cohomology classes 
$\{x^{p_1}_{1}\dots x^{p_n}_{n}\}$ and $\{x^{l_1}_{1}\dots x^{l_n}_{n}\}$ such that 
$$x^{i_1}_{1}\dots x^{i_n}_{n}=\big(x^{p_1}_{1}\dots x^{p_n}_{n}\big)\big(x^{l_1}_{1}\dots x^{l_n}_{n}\big),\,\,\,\,\,\,r_p=\sum_{i}^{n}p_i\vert x_{i}\vert -1,\,\,\,\, r_q=\sum_{i}^{n}l_i\vert x_{i}\vert -1$$
$$p_1\geq 0, \dots,p_n\geq 0,\,\,\,\,\,\,\,\,\,\,\,\,\,\,\,\,\,\, l_1\geq 0, \dots,l_n\geq 0.$$

\medskip
It well-known that if $M$ is the formal dimension of the $F_0$-space $X$, then, thanks to the Poincar\'{e} duality (\cite{FHT}, \S 38), we have an isomorphism of vector spaces 
$$\phi:W_{r_s}\to W_{M-2-r_s}.$$
So if $\big\{w_{(i,r_s)}\big\}_{1\leq i\leq n_s}$ is a basis for $W_{r_s}$, then $\big\{\phi(w_{(i,r_s)})=w^*_{(i,r_s)}\big\}_{1\leq i\leq n_s}$ is a basis  for $W_{M-2-r_s}$,  called the dual basis. Consequently,  we can choose a basis for  $W_{}$ on the form 
\begin{equation}\label{x4}
\B=\Big\{w_{(1,r_s)},\dots,w_{(n_s,r_s)};w^*_{(1,r_s)},\dots,w^*_{(r_{s},r_s)},\mu\Big\}_{r_1\leq r_s\leq \frac{M-2}{2}},
\end{equation} 
where $W_{M-1}=\langle \mu\rangle$. Moreover, due to (Theorem 2, \cite{S}), we have
\begin{equation}\label{x5}
\delta(\mu)=\frac{1}{2}\sum_{r_s,t}^{} [w_{(t,r_s)},w^*_{(t,r_s)}],\,\,\,\,\,\,\,\,\,\,\,\,\,\,\,1\leq s\leq m ,\,\,\,\,\,\, \,\,\,\,\,\,\,\,\,\,1\leq t\leq n_{r_s}. 
\end{equation}
Note that the integer  $M$ is even (see Remark \ref{r4}), and if $r_p<r_q$, then $\vert w^*_{t,r_q}\vert< \vert w^*_{t,r_p}\vert.$ 

\medskip
The following  result plays a crucial role afterwards.
\begin{lemma}\label{ll1}
	Let $(\L(W_{\rm{}}),\delta)$ be  the Quillen model of an $F_0$-space $X$ of formal dimension $M$.	For every  $w^*_{(j,r_s)}\in\B$,  there exists     $w^*_{(k,r_{\sigma})}\in \B$ such that 
\begin{equation}\label{x42}
\delta(w^*_{(k,r_{\sigma})})=\beta_{(k,r_{\sigma})}[w_{(s_1,r_p)},w^*_{(j,r_s)}]+\Theta_{(k,r_{\sigma})},
\end{equation}
where   $\Theta_{(k,r_{\sigma})}$ is a linear combination of $2$-brackets where  $w_{(s_1,r_p)}$ and $w^*_{(j,r_s)}$ are not involved. Moreover, $w_{(s_1,r_p)}$ is a cycle.
\end{lemma} 
\begin{proof}
	First,  recall that  $\vert w_{(j,r_s)}\vert=r_s$ and $\vert w^*_{(j,r_s)}\vert=M-2-r_s$. 
	Next,  by (\ref{x12}) and (\ref{y2}) we know that  to $w_{(j,r_s)}$ and  $w^*_{(j,r_s)}$ correspond two non-trivial  cohomology classes   
	$\{x^{t_1}_{s_1}\dots x^{t_h}_{s_h}\}$ and $\{x^{i_1}_{j_1}\dots x^{i_k}_{j_k}\}$  in the Sullivan model $(\Lambda V,\partial)$,   such that 
	$$\vert x^{t_1}_{s_1}\dots x^{t_h}_{s_h}\vert=\vert w_{(j,r_s)}\vert+1=r_s+1,\,\,\,\,\,\,\,\,\,\,\,\,\,\,\,\,\,\,\,\,\,\,\,\,\,\,\,\,\,\,\,\,\,\,\,\,\,\,\,\,\,\,\,\,\,\,\,\vert x_{s_1}\vert\leq  \dots\leq \vert x_{s_h}\vert.$$
	$$\vert x^{i_1}_{j_1}\dots x^{i_k}_{j_k}\vert=\vert w^*_{(j,r_s)}\vert+1=M-1-r_s,\,\,\,\,\,\,\,\,\,\,\,\,\,\,\,\,\,\,\,\,\,\,\,\,\,\,\,\,\,\,\,\,\,\,\,\,\vert x_{j_1}\vert\leq  \dots\leq \vert x_{j_k}\vert.$$
	Here  we can assume $t_1\geq 1,  \dots, t_h\geq 1$ and $i_1\geq 1,  \dots,i_k\geq 1$. Note that if the generator  $w_{(j,r_s)}$ is a cycle, then  the corresponding element in $(\Lambda V,\partial)$ is the cohomology class   $\{x^{}_{s_1}\}$.
	
Next,  Poincar\'{e} duality implies that the 
	multiplication
	$$H^{r_s+1}(\Lambda V)\times H^{M-1-r_s}(\Lambda V)\to H^{M}(\Lambda V),$$
	sending 
$	\Big(\{x^{t_1}_{s_1}\dots x^{t_h}_{s_h}\}; \{x^{i_1}_{j_1}\dots x^{i_k}_{j_k}\}\Big)$ to $\{x^{t_1}_{s_1}\dots x^{t_h}_{s_h}.x^{i_1}_{j_1}\dots x^{i_k}_{j_k}\},$
	is non-degenerate.  It follows that   $x_{s_i}(x^{i_1}_{j_1}\dots x^{i_k}_{j_k})$ is not a coboundary for every $1\leq i\leq h$. As a result, we must have  a generator  $w^*_{(k,r_\sigma)}$ corresponding to cohomology class $\{x_{s_i}(x^{i_1}_{j_1}\dots x^{i_k}_{j_k})\}$ such that $\delta(w^*_{(k,r_{\sigma})})$ satisfies the following formula 
$$	\delta(w^*_{(k,r_{\sigma})})=\beta_{(k,r_{\sigma})}[w_{(s_1,r_p)},w^*_{(j,r_s)}]+\Theta_{(k,r_{\sigma})},$$
where $w_{(s_1,r_p)}$ corresponds to $x_{s_1}$ which implies that $w_{(s_1,r_p)}$ is a cycle.

\noindent Finally, from the formula (\ref{y2}), it is clear that  $\Theta_{(k,r_{\sigma})}$ is a linear combination of $2$-brackets where  $w_{(s_1,r_p)}$ and $w^*_{(j,r_s)}$ are not involved.
\end{proof}
\begin{remark}\label{rr4} 
In the cohomology class  $\{x^{t_1}_{s_1}\dots x^{t_h}_{s_h}\}$   corresponding  to $w_{(j,r_s)}$,	we  might have $$\vert x_{s_1}\vert=  \dots= \vert x_{s_\tau}\vert,\,\,\,\,\,\,\,\,\,\,\,\,\,\,\,1\leq \tau\leq h.$$ 
In this case,  the formula (\ref{x42}) can be written as follows
	\begin{equation*}\label{x44}
		\delta(w^*_{(k,r_{\sigma})})=\beta_{1}[w_{(s_1,r_p)},w^*_{(j,r_s)}]+\sum_{ j'\neq j\,,\,i>1}^{h} \beta_{i}[w_{(s_i,r_p)},w^*_{(j',r_s)}]+\Theta_{(k,r_{\sigma})},
	\end{equation*}
furthermore, we have the following facts.
\begin{enumerate}
		\item Since $\vert x_{s_1}\vert\leq  \dots\leq \vert x_{s_h}\vert$, we deduce that $\Theta_{(k,r_{\sigma})}$ is a linear combination of $2$-brackets of the form $[w_{(a,b)},w_{(c,d)}]$ such that  $$r_p<\vert w_{(a,b)}\vert \leq \vert w_{(c,d)}\vert <M-2-r_s.$$
	\item All the generators  $w^*_{(j,r_s)}$ and $w^*_{(j',r_s)}$, where $j'\neq j$, are distinct and have the same degree $M-2-r_s$.
	\item All the generators  $w_{(s_i,r_p)}, 1\leq i\leq h$, are distinct cycles with  $\vert w_{(s_i,r_p)}\vert=r_p$.
	\item All the rationals $\beta_{i}$ are not zero.
\end{enumerate}
\end{remark}
\begin{remark}\label{rr3}
 A special case of Lemma \ref{ll1} is when $r_s= \frac{M-2}{2}$.  In this case the lemma    still valid for any  generator  $w_{(j,r_s)}$	such that $\delta(w^*_{(j,r_s)})\neq 0$ because the dual of $w^*_{(j,r_s)}$, namely $(w^*_{(j,r_s)})^*$, is $w_{(j,r_s)}.$ 
\end{remark}
\section{Main result}
In all this section,  let $X$ denote an $F_0$-space of formal dimension $M$,  $(\Lambda V,\partial)$  its Sullivan model,  $(\L(W_{\rm{}}),\delta)$   its  Quillen model and  $\B$  the basis of $W$ given in (\ref{x4}). Recall that by  Corollary (\ref{c3}) there exists a rational $a\neq 0$ such that $\alpha(\mu)=a\mu$, where $W_{M-1}=\langle \mu\rangle$.

\medskip
Subsequently,  we  prove   some important lemmas concerning the properties of    $(\L(W_{\rm{}}),\delta)$  needed to establish  the main result in this paper. Indeed, if   $[\alpha]\in \E_{\#}(\L(W))$, then by considering the basis (\ref{x4}) and Remark \ref{rr4}, we can summarize the next steps as follows.
\begin{itemize}
	\item  In Lemma \ref{ll2}, we  show that  $\alpha(w^*_{j,r_s})=aw^*_{j,r_s}$ for all $j$ and $r_s< \frac{M-2}{2}.$
	\item In Lemma \ref{ll3}, we show that $\alpha(w_{j,r_s})=w_{j,r_s}$ for every $j$ and $r_s< \frac{M-2}{2}$.
	\item  Lemmas \ref{ll4} and  \ref{c01}, show that $\alpha(w_{j,\xi})=aw_{j,\xi}$ and $\alpha(w^*_{j,\xi})=aw^*_{j,\xi}$ for every $j$, where $\xi=\frac{M-2}{2}$.
	\item  In  Proposition  \ref{p1}, we show  that $a=1$. 
\end{itemize}
\begin{lemma}\label{ll2}
	 If $[\alpha]\in \E_{\#}(\L(W))$, then   for  $w^*_{j,r_s}\in\B$ such that  $r_s< \frac{M-2}{2}$,   we have $\alpha(w^*_{j,r_s})=aw^*_{j,r_s}$. 
\end{lemma} 
\begin{proof} 
	First, let us  prove that for every $j$ we have
	 \begin{equation}\label{x31}
		\alpha(w^*_{(j,r_1)})=aw^*_{(j,r_1)}.
	\end{equation}
	Indeed,  let $ w^*_{(j,r_1)}\in W_{M-2-r_1}=\big\langle w^*_{(1,r_1)},w^*_{(2,r_1)},\dots \big\rangle$.  Using  Corollary \ref{c3}, we can write 
	\begin{equation}\label{x7}
		\alpha(w^*_{(j,r_1)})=\sum_{i\geq 1}^{}\lambda_{i}w^*_{(i,r_1)},\,\,\,\,\,\,\,\,\,\,\,\,\,\,\,\,\,\,\,\,\,\,\,\,\,\,\,\,\,\,\,\,\,\,\,\lambda_{i}\in\Q.
	\end{equation}
	The formula (\ref{x5})  can be written as 
	\begin{equation}\label{x29}
		\delta(\mu)=\frac{1}{2}\,\,\sum_{ j} [w_{(j,r_1)},w^*_{(j,r_1)}]+\frac{1}{2}\sum_{r_1\neq r_s} [w_{(t,r_s)},w^*_{(t,r_s)}].
	\end{equation}	
	Next,  $ w_{(j,r_1)}$ is obviously a cycle as $ w_{(j,r_1)}\in W_{r_1}$.  By Corollary \ref{c4}, it follows that
	\begin{equation}\label{x6}
		\alpha(w_{(j,r_1)})=w_{(j,r_1)}.
	\end{equation}
	On the one hand, 	using (\ref{x7}), (\ref{x29}) and  (\ref{x6}), we get 
	\begin{equation*}\label{x30}
		\alpha(\delta(\mu))=\frac{1}{2}\,\,\sum_{ j}\sum_{i\geq 1}^{}\lambda_{i} [w_{(j,r_1)},w^*_{(i,r_1)}]+\frac{1}{2}\sum_{r_1\neq r_s} [\alpha(w_{(t,r_s)}),\alpha(w^*_{(t,r_s)})].
	\end{equation*}
	On the other hand, by the relation (\ref{x29})  and Corollary \ref{c3}, we have
	$$\delta(\alpha(\mu))=a\delta(\mu)=\frac{a}{2}\,\,\sum_{ j} [w_{(j,r_1)},w^*_{(j,r_1)}]+\frac{a}{2}\sum_{r_1\neq r_s} [w_{(t,r_s)},w^*_{(t,r_s)}].$$
	Since $\alpha(\delta(\mu))=\delta(\alpha(\mu))$, and taking into account that $r_1\neq r_s$, which means that the generator $w^*_{(j,r_1)}$ cannot appear in the expression $\sum_{r_1\neq r_s} [\alpha(w_{(t,r_s)}),\alpha(w^*_{(t,r_s)})]$, it follows that all the rationals $\lambda_{i}$ in (\ref{x7}) are zero except   $\lambda_{1}=a$ showing (\ref{x31}).
	\medskip

 Next, assume by induction that 
 \begin{equation}\label{x15}
 \alpha(w^*_{(j,r_q)})=aw^*_{(j,r_q)},
 \end{equation}
  for all the generators  $w^*_{j,r_q}$ such that $r_{q}<r_{s}$.  Let us prove it for   every generator 
 $$ w^*_{(j,r_s)}\in W_{M-2-r_s}=\big\langle w^*_{(1,r_s)},w^*_{(2,r_s)},\dots\big\rangle.$$ 
 For this purpose, write
 \begin{equation}\label{x40}
 	\alpha(w^*_{(j,r_s)})=\sum_{\tau\geq 1}^{}\lambda_{\tau}w^*_{(\tau,r_s)},\,\,\,\,\,\,\,\,\,\,\,\,\,\,\,\,\,\,\,\,\,\,\,\,\,\,\,\,\,\,\,\,\,\,\,\,\,\,\,\,\,\,\,\,\lambda_{\tau}\in\Q.
 \end{equation}
 By virtue of Lemma \ref{ll1} and Remark \ref{rr4}, there exists $w^*_{(k,r_{\sigma})}$ such that  
 \begin{equation}\label{x39}
 		\delta(w^*_{(k,r_{\sigma})})=\beta_{1}[w_{(s_1,r_p)},w^*_{(j,r_s)}]+\hspace{-2mm}\sum_{ j'\neq j\,,\,i>1}^{h} \hspace{-3mm}\beta_{i}[w_{(s_i,r_p)},w^*_{(j',r_s)}]+\Theta_{(k,r_{\sigma})},
 \end{equation}
where $\beta_1\neq 0$.  As a result,  we obtain 
 \begin{eqnarray}\label{x10}
 	\alpha(\delta(w^*_{(k,r_{\sigma})}))\hspace{-2mm}&=&\hspace{-2mm}\beta_{1}[\alpha(w_{(s_1,r_p)}),\alpha(w^*_{(j,r_s)})]+\hspace{-3mm}\sum^{h}_{ j'\neq j\,,\,i>1}\hspace{-3mm} \beta_{i}[\alpha(w_{(s_i,r_p)}),\alpha(w^*_{(j',r_s)})]+\alpha(\Theta_{(k,r_{\sigma})})\nonumber\\
 	\hspace{-2mm}&=&\hspace{-2mm}\sum_{\tau\geq 1}^{}\lambda_{\tau}\beta_{1}[w_{(s_1,r_p)},w^*_{(\tau,r_s)}]+\hspace{-3mm}\sum^{h}_{ j'\neq j\,,\,i>1} \hspace{-3mm}\beta_{i}[w_{(s_i,r_p)},\alpha(w^*_{(j',r_s)})]+\alpha(\Theta_{(k,r_{\sigma})}).\nonumber
 \end{eqnarray}
Note that, according to Remark \ref{rr4},  all the generators  $w_{(s_i,r_p)}$  are cycles implying that 
$\alpha(w_{(s_i,r_p)})=w_{(s_i,r_p)}$ due to Corollary \ref{c4}.

\noindent Next, as $\vert w^*_{(k,r_{\sigma})}\vert>\vert w^*_{(j,r_s)}\vert$ which implies that $r_{q}<r_{s}$,  using (\ref{x15}) and (\ref{x39}), we get
 $$\delta(\alpha(w^*_{(k,r_{\sigma})}))=a\delta(w^*_{(k,r_{\sigma})})=a\beta_{1}[w_{(s_1,r_p)},w^*_{(j,r_s)}]+\hspace{-2mm}\sum_{ j'\neq j\,,\,i>1}^{h} \hspace{-3mm}a\beta_{i}[w_{(s_i,r_p)},w^*_{(j',r_s)}]+a\Theta_{(k,r_{\sigma})}.$$
 Since $\alpha(\delta(w^*_{(i,r_q)}))=\delta(\alpha(w^*_{(i,r_q)}))$ and taking into account  that the bracket $[w_{(s_1,r_p)},w^*_{(j,r_s)}]$  does not appear in the expressions  (see Remark \ref{rr4})
 $$\sum_{ j'\neq j\,,\,i>1}^{h} \hspace{-3mm}a\beta_{i}[w_{(s_i,r_p)},w^*_{(j',r_s)}]\,\,\,\,\,\,\,\,\text{ and }\,\,\,\,\,\,\,\,\alpha(\Theta_{(k,r_{\sigma})}),$$
  we deduce that all the coefficients 
 $\lambda_{\tau}$ in (\ref{x40}) are nil  except  $\lambda_{1}\beta_{1}=a\beta_{1}$ and because $\beta_{1}\neq 0$, we obtain  $\alpha(w^*_{(j,r_s)})=aw^*_{(j,r_s)}$. 
\end{proof}
	\begin{lemma}\label{ll3}
	 If $[\alpha]\in \E_{\#}(\L(W))$, then   for every  $w_{j,r_s}\in\B$, where $r_s< \frac{M-2}{2}$,   we have $\alpha(w_{j,r_s})=w_{j,r_s}.$
\end{lemma} 
\begin{proof}
First, we know from   Corollary \ref{c4} that if $\delta(w_{j,r_s})=0$, then $\alpha(w_{j,r_s})=w_{j,r_s}$, therefore we can suppose that $\delta(w_{j,r_s})\neq0$. Secondly,   recall that from the formula (\ref{x5}) we can write 
\begin{equation*}\label{x18}
	\delta(\mu)=\frac{1}{2}[w_{(j,r_s)},w^*_{(j,r_s)}]+\frac{1}{2}\,\,\sum_{ t\neq j}[w_{(t,r_s)},w^*_{(t,r_s)}]. 
\end{equation*}
As a result, we get
\begin{equation*}\label{x24}
	\alpha(\delta(\mu))=\frac{1}{2}[\alpha(w_{(j,r_s)}),\alpha(w^*_{(j,r_s)})]+\frac{1}{2}\,\,\sum_{ t\neq j}[\alpha(w_{(t,r_s)}),\alpha(w^*_{(t,r_s)})],
\end{equation*}
and because $r_s< \frac{M-2}{2}$,  Lemma \ref{ll2} implies that 
$$\alpha(w^*_{(j,r_s)})=aw^*_{(j,r_s)},\,\,\,\,\,\,\,\alpha(w^*_{(t,r_s)})=aw^*_{(t,r_s)},\,\,\,\,\,\,\,\forall t\neq j.$$ 
Next,  by Corollary \ref{c3},   we  can  write 
\begin{equation*}\label{x36}
\alpha(w_{(j,r_s)})=\sum_{i}^{}\rho_{i}w_{(i,r_s)},\,\,\,\,\,\,\,\,\,\,\,\,\,\,\,\,\,\,\,\,\,\rho_{i}\in\Q,
\end{equation*}
implying that  
\begin{equation}\label{x16}
	\alpha(\delta(\mu))=\frac{a}{2}\rho_{j}[w_{(j,r_s)},w^*_{(j,r_s)}]+\frac{a}{2}\sum_{i\neq j}^{}\rho_{i}[w_{(i,r_s)},w^*_{(j,r_s)}]+\frac{a}{2}\,\,\sum_{ t\neq j}[\alpha(w_{(t,r_s)}),w^*_{(t,r_s)}].
\end{equation}
	Finally,  using Corollary \ref{c4}, we obtain
	\begin{equation}\label{x17}
	\delta(\alpha(\mu))=	a\delta(\mu)=\frac{a}{2}[w_{(j,r_s)},w^*_{(j,r_s)}]+\frac{a}{2}\,\,\sum_{ t\neq j}[w_{(t,r_s)},w^*_{(t,r_s)}]. 
	\end{equation}
	Since $\alpha(\delta(\mu))=\delta(\alpha(\mu))$ and  $w^*_{(j,r_s)}\neq w^*_{(t,r_s)}$,  comparing (\ref{x16}) and (\ref{x17}), it follows that 
	$\rho_{i}=0$ for all $i\neq j$  and  $\rho_{j}=1$. 
	Hence, $\alpha(w_{(j,r_s)})=w_{(j,r_s)}.$
\end{proof}
\begin{lemma}\label{ll4}
If $[\alpha]\in \E_{\#}(\L(W))$, then   $\alpha(w^*_{(j,\xi)})=aw^*_{(j,\xi)}$, where  $\xi= \frac{M-2}{2}$.
\end{lemma} 
\begin{proof} 
 By virtue of Lemma \ref{ll1} and Remark \ref{rr4}, there exists $w^*_{(k,r_\sigma)}$ such that  
\begin{equation}\label{x21}
		\delta(w^*_{(k,r_{\sigma})})=\beta_{1}[w_{(s_1,r_p)},w^*_{(j,\xi)}]+\hspace{-4mm}\sum_{ j'\neq j\,,\,i>1}^{h} \hspace{-3mm}\beta_{i}[w_{(s_i,r_p)},w^*_{(j',\xi)}]+\Theta_{(k,r_{\sigma})},
\end{equation}
where the generators $w_{(s_i,r_p)}$ are cycles implying that  $\alpha(w_{(s_i,r_p)})=w_{(s_i,r_p)}$ for all $1\leq i\leq h$. Next,  since  that a basis of $W_{\xi}$ is formed by the generators $w_{(i,\xi)}$ and their duals   $w^*_{(i,\xi)}$ because in this case we have $\vert w_{(i,\xi)}\vert=\vert w^*_{(i,\xi)}\vert=\xi=\frac{M-2}{2}$, by Corollary \ref{c4},   we  can  write 
\begin{equation}\label{x25}
	\alpha(w^*_{(j,\xi)})=\sum_{i\geq 1}^{}\mu_{i}w^*_{(i,\xi)}+\sum_{\tau \geq 1}^{}\gamma_{\tau}w_{(\tau,\xi)},
\end{equation}
As a result,  we obtain 
\begin{eqnarray}\label{xx27}
		\alpha(\delta(w^*_{(k,r_{\sigma})}))&=&\sum_{i}^{}\mu_{i}\beta_{1}[w_{(s_1,r_p)},w^*_{(i,\xi)}]+\sum_{\tau}^{}\gamma_{\tau}\beta_{1}[w_{(s_1,r_p)},w_{(\tau,\xi)}]\nonumber\\
		&+&\hspace{-4mm}\sum_{ j'\neq j\,,\,i>1}^{h} \hspace{-3mm}\beta_{i}[w_{(s_i,r_p)},\alpha(w^*_{(j',\xi)})]+\alpha(\Theta_{(k,r_{\sigma})})\nonumber.
\end{eqnarray}
Next, as $\vert w^*_{(k,r_\sigma)}\vert>\vert w^*_{(j,\xi)}\vert=\xi=\frac{M-2}{2}$, it follows that $r_\sigma< \xi$.  Thus,  using Lemma \ref{ll2} and the relation (\ref{x21}), we get
$$\delta(\alpha(w^*_{(i,r_q)}))=a\delta(w^*_{(i,r_q)})=a\beta_{1}[w_{(s_1,r_p)},w^*_{(j,\xi)}]+\hspace{-4mm}\sum_{ j'\neq j\,,\,i>1}^{h} \hspace{-3mm}a\beta_{i}[w_{(s_i,r_p)},w^*_{(j',\xi)}]+a\Theta_{(k,r_{\sigma})}.$$
Since $\alpha(\delta(w^*_{(i,r_q)}))=\delta(\alpha(w^*_{(i,r_q)}))$ and taking into account that the bracket  $[w_{(s_1,r_p)},w^*_{(j,\xi)}]$  does not appear in the expression $\alpha(\Theta_{(i,r_q)})$, according to Remark \ref{rr4}, 
we deduce that all the coefficients 
$\mu_{i}$ and $\gamma_{\tau}$ in (\ref{x25}) are nil  except  $\mu_{j}\beta_{1}=a\beta_{1}$ implying that $\mu_{j}=a$ because $\beta_{1}\neq 0$. Hence, $\alpha(w^*_{(j,\xi)})=aw^*_{(j,\xi)}$  
\end{proof}
\begin{lemma}\label{c01}
	 If $[\alpha]\in \E_{\#}(\L(W))$, then   for every  $w_{(j,\xi)}\in\B$,   we have $\alpha(w_{(j,\xi)})=aw_{(j,\xi)}$.
\end{lemma} 
\begin{proof} 
	The proof is as in  Lemma \ref{ll4} after taking into consideration Remark \ref{rr3}.
\end{proof}
\begin{proposition}\label{p1}
	If  $(\L(W_{\rm{}}),\delta)$ is  the Quillen model of an $F_0$-space of formal dimension $M$,  then the group  $\E_{\#}(\L(W))$is trivial.
\end{proposition} 
\begin{proof} 
	It suffices to prove that the rational $a$ given  in Lemmas  \ref{ll2},  \ref{ll4} and \ref{c01} satisfies  $a=1$. Indeed, first the formula (\ref{x5})  can be written as 
	\begin{equation*}\label{x32}
		\delta(\mu)=\frac{1}{2}\sum_{j}[w_{(j,\xi)},w^*_{(j,\xi)}]+\frac{1}{2}\sum_{r_p< \xi}\,\,\sum_{t}[w_{(t,r_p)},w^*_{(t,r_p)}]. 
	\end{equation*}	
	It follows that 
	\begin{equation*}\label{x22}
		\alpha(\delta(\mu))=\frac{1}{2}\sum_{j}[\alpha(w_{(j,\xi)}),\alpha(w^*_{(j,\xi)})]+\frac{1}{2}\sum_{r_p< \xi}\,\,\sum_{t}[\alpha(w_{(t,r_p)}),\alpha(w^*_{(t,r_p)})].
	\end{equation*}
Now, for all $t$ and $r_p< \xi$, Lemmas \ref{ll2} and \ref{ll3} yield the following
	$$\alpha(w_{(t,r_p)})=w_{(t,r_p)},\,\,\,\,\,\,\,\,\,\,\,\,\,\,\,\,\,\,\,\,\,\,\,\,\,\,\,\,\,\,\,\,\,\,\,\alpha(w^*_{(t,r_p)})=aw^*_{(t,r_p)},$$
and for for all $t$, by lemmas \ref{ll4} and Corollary \ref{c01},   we have
$$\alpha(w_{(t,\xi)})=aw_{(t,\xi)},\,\,\,\,\,\,\,\,\,\,\,\,\,\,\,\,\,\,\,\,\,\,\,\,\,\,\,\,\,\,\,\,\,\alpha(w^*_{(t,\xi)})=aw^*_{(t,\xi)},$$
Therefore, on the one hand, we have
	\begin{equation*}\label{x35}
		\alpha(\delta(\mu))=\frac{1}{2}\sum_{j}a^2[w_{(i,\xi)},w^*_{(j,\xi)}]+\frac{a}{2}\sum_{r_p< \xi}\,\,\sum_{t}[w_{(t,r_p)},w^*_{(t,r_p)}].
	\end{equation*}
	On the other hand, by the relation (\ref{x32}) and Corollary \ref{c3} we have
	$$\delta(\alpha(\mu))=\frac{a}{2}\sum_{j}[w_{(j,\xi)},w^*_{(j,\xi)}]+\frac{a}{2}\sum_{r_p< \xi}\,\,\sum_{t}[w_{(t,r_p)},w^*_{(t,r_p)}].  $$
	Since $\alpha(\delta(\mu))=\delta(\alpha(\mu))$, it follows that $a^2=a$ and as $a\neq 0$, it follows that $a=1.$
\end{proof} 
Now we are able to announce the main result in this paper.
\begin{theorem}\label{t1}
	If $X$ is   an $F_0$-space,  the $\E_{\#}(X)$ is finite.
\end{theorem} 
\begin{proof}
 It suffices to apply Proposition \ref{p1} and the identification (\ref{f13}).
\end{proof}

\bibliographystyle{amsplain}

\begin{thebibliography}{10}
	\bibitem{ak1} Arkowitz, M., {\em Problems on Self-homotopy equivalences}, Contemporary Math., 274, (2001), 309-315.
	
	
\bibitem{B3} M. Benkhalifa, \emph{The effect of  cell-attachment on the group of self-equivalences of an elliptic  space}, Michigan Mathematical Journal, Vol. 71(2), 611-617, 2022.	
	
	
	
	
\bibitem{B1}M. Benkhalifa, \emph{On the group of self-homotopy equivalences of an elliptic space},
Proceedings of the American Mathematical Society. Vol.148 (6),2695-2706, 2020
	
	
	
	
	
	
	
	
	
	
	
	
	\bibitem{FHT} Y. F{\'e}lix, S.  Halperin, and J.-C. Thomas, \emph{Rational
		homotopy theory}, Graduate Texts in Mathematics, Vol. 205, Springer-Verlag,
	New York, 2001.
	
	\bibitem{H}	Halperin S., {\em Finiteness in the minimal models of Sullivan, } Trans. Amer. Math. Soc. 230, 269-331, 1977.
	
	
	
	\bibitem{G}	G. Lupton, {\em A note on the conjecture of S. Halperin}, Lectures Notes
	in Mathematics, vol. 1440, Springer-Verlag, 148-163, 1990.
	
	\bibitem{Mar}  K. Maruyama,
	{\em Localization of a certain subgroup of self-homotopy equivalences,}
	Pacific J. Math. {136},  293-301, 1989.
	
\bibitem{M} W. Meier,	{\em   Rational universal fibrations and flag manifolds}. Math. Ann., 258(3):329–340,
	1981.
	
	
	
	\bibitem{N} J.  Neisendorfer  and  T.  Miller,
	{\em Formal and coformal spaces,}
	Illinois J. Math. Vol. 205 (4), 565-580, 1978
	

	
\bibitem{S}	J. Stasheff,	{\em   Rational Poincar\'{e} duality spaces }. 	Illinois J. Math. Vol. 27 (1), 104–109,
1983.	

\end{thebibliography}

\end{document}